\documentclass{article}
\usepackage{amsmath}
\usepackage{amsfonts}
\usepackage{amssymb}
\usepackage{euscript}
\input epsf

\newtheorem{thm}{Theorem}[section]

\newtheorem{prop}{Proposition}[section]

\def\dj{d\kern-0.4em\char"16\kern-0.1em}

\def\oph{\oplus^h}
\def\omh{\ominus^h}

\title{\bf The deformed exponential \\functions of two variables
}

\author{
{\bf Miomir S. Stankovi\'c}\\
Dep. of Math., Faculty of Occupational Safety\\
{\bf Sladjana D. Marinkovi\'c}\\
Dep. of Math., Faculty of Electronic Engineering\\
{\bf Predrag M. Rajkovi\'c}\\
Dep. of Math., Faculty of Mechanical Engineering\\
[2mm] {\bf   University of Ni\v s,\ Serbia }
}

\date{}

\begin{document}
\maketitle

\medskip

\noindent{\bf Abstract.} In the recent development in a various disciplines of physics, it is noted the need for including the deformed
versions of the exponential functions. In this paper, we consider the deformations which have two purposes: to have them like special
cases, and, even more, to acquit their inauguration from mathematical point of view. Really, we will show interesting differential and
difference properties of our deformations which are important in forming and explanation of continuous and discrete models of numerous
phenomena.
\\[1mm] {\bf Keywords:} Exponential function, differential operator, difference operator.
\\[2mm]

\section{Introduction}

According as solving of the concrete real problems was leading to the definition of exponential function at last quarter of seventeenth
century, new circumstances and challenges at the of twentieth century  required its generalizations and deformations.
One-parameter deformations of exponential function have been proposed in the context of non-extensive statistic mechanics 
(see \cite{Tsallis1}, \cite{Tsallis2}, \cite{Nivanen}, and \cite{Borges}) relativistic statistical theory 
(see, for example, \cite{Kaniadakis1} and \cite{Kaniadakis2}) and quantum--group theory {\cite{Abe}.} 

C. Tsallis introduced his analogy of exponential function \cite{Tsallis1} 
in 1988.,  Afterwards, considering a connection between the generalized entropy and theory of quantum groups, S. Abe has defined in 1997.
his deformation, and, at last, G. Kaniadakis \cite{Kaniadakis1} 
proposed in 2001. a new one-parameter deformation for the exponential function.

One-parameter deformations of the exponential functions attracted attention of researchers from various scientific areas because of its
successful role in description of fractal structured systems, non-regular diffusion, thermodynamical and gravitational like systems,
optimization algorithm, statistical conclusions and theory of probability, several complex systems, and so on (see \cite{Tsallis2}). 



A lot of natural phenomena have both discrete and continual aspects. The classical mathematics has difficulties in expressing discrete and
continuous structure in one model. That's why the relations between them establish in two directions: (1) from discrete to continual in the
sense of limits (for example, from differences to derivative by limits); (2) from continual to discrete, like approximation of continual
phenomena (for example, the continuous function and its expansion by Newton series).

The areas where deformations of the exponential functions have been treated basically along three (complementary) directions:
formal mathematical developments \cite{Tsallis1}, \cite{Tsallis2}, 
\cite{Kaniadakis1}, \cite{Kaniadakis2}, 
and \cite{Abe1}; 
observation of consistent concordance with experimental (or natural) behavior \cite{Kaniadakis1}, 
and theoretical physical developments \cite{Abe}. 
At last, we will mention that a very interesting discussion about adequacy of introducing the generalizations of known functions can be found in \cite{Qadir}.


In this paper, using formal mathematical approach, we introduce two variants of deformed  exponential function of two variables
to express discrete and continual behavior by one model. It has such differential i difference properties which able us to do it.
In this function someone can recognize well-known generalizations and deformations like the special cases.

The paper is organized as follows: after sections devoted to introduction and preliminaries, we establish the deformed exponential
function of two variables in the third section.
In the last two sections, we examine its difference and differential properties.
We will prove that this functions appear as the eigenfunctions of the difference and differential operators associated with eigenvalues
by the first and second variable respectively.

\section{Preliminaries: powers and differences}

Let $h\in\mathbb R\setminus\{0\}$. The generalized integer powers of real numbers have important role in modern theoretical considerations
(see, for example, \cite{Kac}, \cite{Riordan}). 
In that manner we firstly introduce {\it backward and forward integer power} given by

\begin{eqnarray*}
z^{(0,h)}&=&1,\qquad
z^{(n,h)}
=\prod_{k=0}^{n-1}(z-kh)\quad (n\in\mathbb N),\\
z^{[0,h]}&=&1,\qquad
z^{[n,h]}
=\prod_{k=0}^{n-1}(z+kh)\quad (n\in\mathbb N).\\
\end{eqnarray*}
The {\it central integer power} is given by
\begin{eqnarray*}
z^{\langle 0,h\rangle}&=&1,\nonumber\\
z^{\langle n,h\rangle}&=&
\left\{\begin{array}{ll}
\prod\limits_{k=0}^{m-1}(z-2kh)(z+2kh)&\ (n=2m,\ m\in\mathbb N),\\
z\prod\limits_{k=0}^{m-1}\bigl(z-(2k+1)h\bigr)\bigl(z+(2k+1)h\bigr)&\ (n=2m+1,\ m\in\mathbb N_0).
\end{array}\right.
\end{eqnarray*}

For $n\in\mathbb N$, the connection with the previous defined generalized powers is given by
\begin{eqnarray}
z^{(n,h)}&=&z^{[n,-h]},\qquad z^{[n,h]}\ \ =\ z^{(n,-h)},\label{-h}\\
z^{\langle n,h\rangle}&=&z^{\langle n,-h\rangle},\qquad
z^{\langle n,h\rangle}\ =\ z\bigl(z+(n-2)h\bigr)^{(n-1,\ 2h)}.\label{z<n}
\end{eqnarray}
Also, the following holds:
$$
z^{\langle n,h\rangle}
=\left\{\begin{array}{ll}
z^{(m,2h)}\ z^{[m,2h]}\ &\ (n=2m,\ m\in\mathbb N),\\
z\ (z-h)^{(m,2h)}(z+h)^{[m,2h]}\ &\ (n=2m+1,\ m\in\mathbb N_0)\
\end{array}\right.
$$
$$
z^{\langle 2m,h\rangle}z^{\langle 2m+1,h\rangle}
=z\ z^{(2m,h)}\ z^{[2m,h]}.
$$
Some binomial coefficients can be represented over the generalized powers. So, for $n\in\mathbb N$, the following holds:
\begin{eqnarray}\label{binom}
\binom{z/h}{n}
&=&\frac{z(z-h)\cdots\bigl(z-(n-1)h\bigr)}{h^n\ n!}
=\frac{z^{(n,h)}}{h^n\ n!}\ .
\end{eqnarray}

Consider the $h$--difference operators \cite{Jagerman} 
\begin{eqnarray}
\Delta_{z,h} f(z)&=&\frac{f(z+h) - f(z)}{h}\;,\label{forwOp} \\
\nabla_{z,h} f(z)&=&\frac{f(z) - f(z-h)}{h}\;,\label{backOp}  \\
\delta_{z,h} f(z)&=&\frac{f(z+h) -  f(z-h)}{2h}\;.\label{centOp}
\end{eqnarray}
Notice that
$$
\aligned
\nabla_{z,h}f(z)=\Delta_{z,-h}&f(z)=\Delta_{z,h}f(z-h),\\
\delta_{z,-h} f(z)&=\delta_{z,h} f(z).
\endaligned
$$
We can prove that their acting on integer generalized powers is given by:
\begin{eqnarray}
\Delta_{z,h}\,z^{(n,h)}&=&nz^{(n-1,h)},\label{Delta_z^n}\\
\nabla_{z,h}\,z^{[n,h]}&=&nz^{[n-1,h]},\label{nabla_z^n}\\
\delta_{z,h}\,z^{\langle n,h\rangle}&=&nz^{\langle n-1,h\rangle}. \label{delta_z^n}
\end{eqnarray}

\section{The deformed exponential functions}

Let $h\in\mathbb R\setminus\{0\}$. In this section we present two deformations of exponential
function of two variables depending of parameter $h$.

\bigskip

{\bf 1.} Let us define function $(x,y)\mapsto e_h(x,y)$ by
\begin{equation}\label{eh}
e_h(x,y)=(1+hx)^{y/h} \qquad (x\in \mathbb C\setminus\{-1/h\},\ y\in \mathbb R).
\end{equation}
Since
$$
\lim_{h\to 0}e_h(x,y)=e^{xy}\ ,
$$
we may take it for a one--parameter deformation of exponential function of two
variables.

If $h=1-q$ $(q\ne 1)$ and $y=1$, the function (\ref{eh}) becomes
$$
e_{1-q}(x,1)=\bigl(1+(1-q)x\bigr)^{1/(1-q)}\ ,
$$
i.e.,\ $e_{1-q}(x,1)=e^x_q$, where $e^x_q$ is Tsallis
$q$--exponential function [1] defined by
$$
e^x_q  = \left\{
\begin{array}{cc}
\bigl(1+(1-q)x\bigr)^{1/(1-q)}\ ,\quad &1+(1-q)x>0\ ,\\ \\
0\ , &\textrm{otherwise},
\end{array}
\right. \qquad (x\in\mathbb R).
$$

If $h=p-1$ $(p\ne 1)$ and $x=1$, the function (\ref{eh}) becomes
$$
e_{p-1}(1,y)=p^{y/(p-1)},
$$
i.e. function considered as generalization of the standard exponential function in context of quantum group formalism \cite{Johal}. 

Notice that function (\ref{eh}) can be written in the form
$$
e_h(x,y)=\exp\Bigl(\frac yh\ \ln(1+hx)\Bigr).
$$
Hence, we can define deformation function $x\mapsto \{x\}_h$ by
\begin{equation}
\{x\}_h=\frac 1h\ \ln(1+hx)\qquad (x\in \mathbb C\setminus\{-1/h\})\label{{x}_h}
\end{equation}
and, in the same manner,  $(x, y)\mapsto g_h(x,y)$ by
\begin{equation*}
g_h(x,y)=\{x\}_h\ y=\frac yh\ \ln(1+hx)\qquad (x\in \mathbb C\setminus\{-1/h\},\ y\in \mathbb R).
\end{equation*}
So, the following holds:
\begin{equation}
e_h(x,y)=e^{\{x\}_h \;y},\qquad\text{i.e.}\qquad e_h(x,y)=e^{g_h(x,y)}.\label{deform_eh}
\end{equation}

We can show that function (\ref{eh}) holds on some basic properties of exponential function.

\begin{prop} For $x\in \mathbb C\setminus\{-1/h\}$ and $y\in \mathbb R$ the following holds:
$$
\aligned
e_h(x,y)&>0\qquad (x<-1/h\quad \text{for}\ h<0\quad \text{or}\quad x>-1/h\quad \text{for}\ h>0),\\
e_h(0,y)&=e_h(x,0)=1,\\
e_{-h}(x,y)&=e_h(-x,-y)\qquad (x\ne 1/h),\\
e_h(x,y_1+y_2)&=e_h(x,y_1)e_h(x,y_2).
\endaligned
$$
\end{prop}

Notice that the additional property holds on in regard to the second variable only. However, by treatment of the first variable, the
following holds:
$$
e_h(x_1,y)e_h(x_2,y)=e_h\bigl(x_1+x_2+hx_1x_2,y\bigr)\ .
$$
This equality suggests (see [3], [4]) 
introducing a generalization of the sum operation~\footnote{Notice that generalized summing was already using in analyzing of properties of
functions,  see \cite{aczel},\cite{Pompeiu},\cite{Strauss}} 
\begin{equation}\label{e_hplus}
x_1\oplus_h x_2=x_1+x_2+hx_1x_2\ .
\end{equation}
This operation is commutative, associative and 0 is its neutral. For $x\ne -1/h$ the $\ominus_h$--\ opposite number exists as
$$
\ominus_h x=\frac{-x}{1+hx}
$$
and $x\oplus_h (\ominus_h x)=0$ is valid. Hence, $(I,\oplus_h)$ is abelian group, where \linebreak $I=(-\infty,-1/h)$ for $h<0$ or
$I=(-1/h,+\infty)$ for $h>0$ (see \cite{Benghia}). 
In this way, the $\ominus_h$--subtraction can be defined by
\begin{equation*}\label{e_hminus}
x_1\ominus_h x_2=x_1\oplus_h(\ominus_h x_2)=\frac{x_1-x_2}{1+hx_2}\qquad \Bigl(x_2\ne -\frac 1h\Bigr).
\end{equation*}
With respect to (\ref{{x}_h}) we can prove the next equality
for $x_1,x_2\in I$:
$$
\{x_1\}_h+\{x_2\}_h=\{x_1\oplus_h x_2\}_h.
$$

Now, we will prove a few new properties of function (\ref{eh}).

\begin{thm} For $x_1,x_2\in \mathbb C\setminus\{-1/h\}$ and  $y\in \mathbb R$, the following is valid:
$$
\aligned
e_h(x_1\oplus_h x_2,y)&=e_h(x_1,y)e_h(x_2,y),\\
e_h(x_1\ominus_h x_2,y)&=e_h(x_1,y)e_h(x_2,-y).
\endaligned
$$
\end{thm}

\noindent{\it Proof.} The first equality is following immediately from (\ref{eh}) and (\ref{e_hplus}).
For the second one we can notice the following:
$$
\aligned e_h(\ominus_h x,y)&=e_h\Bigl(\frac{-x}{1+hx},y\Bigr)=\Bigl(1-\frac{x}{1+hx}\Bigr)^{y/h}=\frac{1}{(1+hx)^{y/h}}=e_h(x,-y).\ \Box
\endaligned
$$

\bigskip

{\bf 2.} Let us define function $(x,y)\mapsto \exp_h(x, y)$ by
\begin{equation}\label{exph}
\exp_{h}(x)=\bigl(hx+\sqrt{1+h^2 x^2}\bigr)^{y/h} \qquad (x\in \mathbb C,\ y\in \mathbb R).
\end{equation}
Since
$$
\lim_{h\to 0}\exp_h(x,y)=e^{xy}\ ,
$$
it also can be taken for a one--parameter deformation of exponential function of two
variables.

If $h=\kappa$ and $y=1$, the function (\ref{exph}) becomes
$\kappa$-exponential function
$$
\exp_\kappa(x,1)=\exp_{\{\kappa\}}(x)=\bigl(\sqrt{1+\kappa^2x^2}+\kappa x\bigr)^{1/\kappa},
$$
introduced in \cite{Kaniadakis1}, \cite{Kaniadakis2}. 

The connection between functions (\ref{eh}) and (\ref{exph}) is given by
$$
\exp_{h}(x,y)=e_h\Bigl(x-\frac{1-\sqrt{1+h^2 x^2}}{h}\ ,\ y\Bigr).\\
$$

Having in mind that
\begin{equation}\label{arcsinh}
\text{arcsinh}(hx)=\ln(hx+\sqrt{1+h^2 x^2}),
\end{equation}
notice that (\ref{exph}) can be written in the form
$$
\exp_h(x,y)=\exp\Bigl(\frac yh\ \text{arcsinh}\, hx\Bigr).
$$
In \cite{Kaniadakis2}, 
it was defined the deformation function $x\mapsto \{x\}^h$ by
\begin{equation}
\{x\}^h=\frac 1h\ \text{arcsinh}\, hx \qquad (x\in \mathbb C).\label{{x}^h}
\end{equation}
We can define two--dimensional deformation function  $(x,y)\mapsto g^h(x,y)$ as
\begin{equation*}
g^h(x,y)=\frac yh\ \text{arcsinh}\, hx\qquad (x\in \mathbb C,\ y\in \mathbb R).
\end{equation*}
Now, function (\ref{exph}) can be written as
\begin{equation}
\exp_h(x,y)=e^{\{x\}^hy}\qquad \text{i.e.}\qquad \exp_h(x,y)=e^{g^h(x,y)}.\label{deform_exph}
\end{equation}

We adduce the main properties of the introduced deformed exponential function without proof.

\begin{prop}
For $x\in\mathbb C$ and $y\in \mathbb R$ the following holds:
$$
\aligned
\exp_h(x,y)&>0\qquad (x\in\mathbb R),\\
\exp_h(0,y)&=\exp_h(x,0)=1,\\
\exp_{-h}(x,y)&=\exp_h(x,y),\\
\exp_h(x,y_1+y_2)&=\exp_h(x,y_1)\exp_h(x,y_2).
\endaligned
$$
\end{prop}

As the function (\ref{eh}), this function holds on additional property in regard to the second variable only.
However, according to (\ref{arcsinh}) we have:
$$
\exp_h(x_1,y)\exp_h(x_2,y)
=\exp_h\Bigl(x_1\sqrt{1+h^2 x_2^2}+x_2\sqrt{1+h^2 x_1^2},y\Bigr).
$$
This suggests introducing another generalization of sum operation [5]{
\begin{equation}\label{exp_hplus}
x_1 \oph x_2=x_1\sqrt{1+h^2 x_2^2}+x_2\sqrt{1+h^2 x_1^2}\ .
\end{equation}
Operation $\oph$-- sum is commutative, associative, its neutral is 0 and $\oph$-- opposite for $x\in\mathbb R$ \ is\ $-x$. Thus, $(\mathbb
R,\oph)$ is abelian group, and $\omh$-- subtraction can be defined by
$$
x_1 \omh x_2=x_1 \oph (-x_2)=x_1\sqrt{1+h^2 x_2^2}-x_2\sqrt{1+h^2 x_1^2}\ .
$$
Related to (\ref{{x}^h}), we can prove (see [5]) the next equality:
$$
\{x_1\}^h + \{x_2\}^h=\{x_1\oph x_2\}^h.
$$

With respect of operation $\oph$, function (\ref{exph}) has the following properties:

\begin{thm} For $x_1,x_2\in \mathbb C$ and  $y\in \mathbb R$, the following is valid:
$$
\aligned
\exp_h(x_1\oph x_2,y)&=\exp_h(x_1,y)\exp_h(x_2,y),\\
\exp_h(x_1\omh x_2,y)&=\exp_h(x_1,y)\exp_h(-x_2,y)\\
&=\exp_h(x_1,y)\exp_h(x_2,-y).
\endaligned
$$
\end{thm}

\section{Expansions and difference properties of deformed exponential functions}

In this section we consider the expansions of introduced deformed exponential functions. Related to this expansions, we show that functions
$e_h(x,y)$ and $e_{-h}(x,y)$ are eigenfunctions of operators $\Delta_{y,h}$ and $\nabla_{y,h}$ with eigenvalue $x$. Also, function
$\exp_h(x,y)$ is eigenfunction of operator $\delta_{y,h}$ with eigenvalue $x$.

{
\begin{thm}
For functions $(x,y)\mapsto e_h(x,y)$ and $(x,y)\mapsto e_{-h}(x,y)$, the following representations hold respectively:
\begin{eqnarray}
e_h(x,y)&=&\sum_{n=0}^\infty \frac{1}{n!}x^n y^{(n,h)}\qquad (|hx|<1),\label{eh_expan}\\
e_{-h}(x,y)&=&\sum_{n=0}^\infty \frac{1}{n!}x^n y^{[n,h]}\qquad (|hx|<1).\label{e-h_expan}
\end{eqnarray}
\end{thm}

\noindent{\it Proof.}
With respect to well--known expansion
$$
(1+z)^\alpha=\sum_{n=0}^\infty \binom{\alpha}{n}z^n\quad (|z|<1,\
\alpha\in\mathbb R)\ ,
$$
and relation (\ref{binom}), the following holds:
$$
(1+hx)^{y/h}=\sum_{n=0}^\infty \dbinom{y/h}{n}h^n x^n
=\sum_{n=0}^\infty \frac{y^{(n,h)}}{h^n\ n!}h^n x^n\qquad (|hx|<1),
$$
wherefrom we get the required expansion for $e_h(x,y)$. Using (\ref{-h}), we obtain expansion for $e_{-h}(x,y)$.
$\ \Box$

\begin{thm}
The functions $y\mapsto e_h(x,y)$ and $y\mapsto e_{-h}(x,y)$ are the eigenfunctions of operators $\Delta_{y,h}$ and $\nabla_{y,h}$
respectively, with eigenvalue $x$.
\end{thm}

\noindent{\it Proof.} The statement follows from (\ref{Delta_z^n}) and expansion (\ref{eh_expan}). Hence, function $f(y)=e_h(x,y)$
satisfies the difference equation
$$
\Delta_{y,h}\ f(y)=x\ f(y).
$$
In similar way, using (\ref{nabla_z^n}), we can show that $f(y)=e_{-h}(x,y)$ satisfies the difference equation
$$
\nabla_{y,h}\ f(y)=x\ f(y)\ .\quad\Box
$$

\smallskip

\noindent{\bf Remark.} Above representations can be attained vice versa. Since
$$
\aligned
\Delta_{y,h}e_h(x,y)
&=\frac 1{h}\Bigl((1+hx)^{(y+h)/h}-(1+hx)^{y/h}\Bigr) \\
&=x(1+hx)^{y/h}=xe_h(x,y),
\endaligned
$$
with respect to (\ref{Delta_z^n}), the coefficients in expansion
$$
e_h(x,y)=\sum_{n=0}^\infty \frac{a_n(y,h)}{n!}x^n
$$
have to be
$$
a_n(y,h)=y^{(n,h)}\qquad(n\in\mathbb N_0).
$$

\begin{thm}
The function $y\mapsto \exp_h(x,y)$ is the eigenfunction of operator $\delta_{y,h}$ with eigenvalue $x$.
\end{thm}

\noindent{\it Proof.} For function $\exp_h(x,y)$ the following is valid:
$$
\aligned
\delta_{y,h}\exp_h(x,y)
&=\frac 1{2h}\Bigl(\bigl(hx+\sqrt{1+h^2 x^2}\bigr)^{(y+h)/h}-\bigl(hx+\sqrt{1+h^2 x^2}\bigr)^{(y-h)/h}\Bigr)\\
&=\frac{\bigl(hx+\sqrt{1+h^2 x^2}\bigr)^{\frac yh-1}}{2h}\Bigl(\bigl(hx+\sqrt{1+h^2 x^2}\bigr)^2-1\Bigr)\\
&=x\exp_h(x,y)\ .
\endaligned
$$
Therefore, function $f(y)=\exp_h(x,y)$ satisfies difference equation
$$
\delta_{y,h}f(y)=xf(y).\quad \Box
$$

\begin{thm}
The function $(x,y)\mapsto \exp_h(x,y)$ can be represented as
\begin{equation}\label{exph_expan}
\exp_h(x,y)=\sum_{n=0}^\infty \frac{1}{n!}x^n y^{\langle n,h\rangle}\qquad
\end{equation}
\end{thm}

\noindent{\it Proof.} Let us seek expansion of function (\ref{exph}) in the form
$$
\exp_h(x,y)=\sum_{n=0}^\infty \frac{c_n(y,h)}{n!}x^n.
$$
Forasmuch as
$$
\delta_{y,h}\exp_h(x,y)=x\exp_h(x,y),
$$
it follows
$$
\aligned
\delta_{y,h}\sum_{n=0}^\infty \frac{c_n(y,h)}{n!}x^n &=\sum_{n=0}^\infty \frac{\delta_{y,h}c_n(y,h)}{n!}x^n\\
&=\sum_{n=0}^\infty
\frac{c_n(y,h)}{n!}x^{n+1} =\sum_{n=1}^\infty \frac{c_{n-1}(y,h)}{(n-1)!}x^{n}.
\endaligned
$$
Then the coefficients have to be
$$
\aligned
\delta_{y,h}c_0(y)&=0,\\
\delta_{y,h}c_n(y,h)&=nc_{n-1}(y,h),
\endaligned
$$
wherefrom, according to (\ref{delta_z^n}), we yield:
$$
c_n(y,h)=y^{\langle n,h\rangle}\qquad(n\in\mathbb N_0).\quad \Box
$$

\smallskip

\noindent{\bf Remark.} Notice that in expressions (\ref{deform_eh}) and (\ref{deform_exph}) the deformations of variable $x$ appears, but,
in contrary, in expansions (\ref{eh_expan}) and (\ref{exph_expan}), the deformations of powers of $y$ is present.

\section{Differential properties of deformed exponential functions}

In this section we will look up for differential operators which have deformed exponential functions as eigenfunctions.

In \cite{Kaniadakis1}, 
the deformed $h$--differential and $h$--derivative were defined accordingly with operation (\ref{exp_hplus}):
$$
\aligned
d^h z&=\lim_{u\to z}z\omh u,\\
\frac{df(z)}{d^hz}&=\lim_{u\to z}\frac{f(z)-f(u)}{z\omh u}=\sqrt{1+h^2 z^2}\ \frac{df(z)}{dz}\ .
\endaligned
$$

In this sense we can define deformed $h$--differential and $h$--derivative accordingly with operation (\ref{e_hplus}) (see
\cite{Borges}): 
$$
\aligned
d_h z&=\lim_{u\to z}z\ominus_h u,\\
\frac{df(z)}{d_hz}&=\lim_{u\to z}\frac{f(z)-f(u)}{z\ominus_h u}\ .
\endaligned
$$
With respect to (\ref{e_hminus}) we have
$$
\frac{df(z)}{d_hz}=\lim_{u\to z}\frac{f(z)-f(u)}{\dfrac{z-u}{1+hu}}=(1+hz)\dfrac{df(z)}{dx}\ .
$$

\begin{thm}
The function $x\mapsto e_h(x,y)$ is the eigenfunction of operator $\dfrac d{d_h x}$ with eigenvalue $y$.
\end{thm}

\noindent{\it Proof.} Let us apply differential operator $\dfrac{\partial}{\partial x}$ on function $e_h(x,y)$. Firstly, we have
$$
\aligned \frac{\partial}{\partial x}e_h(x,y)&=\sum_{n=1}^\infty \dfrac{1}{(n-1)!}\ x^{n-1} y^{(n,h)}\\
&=\sum_{n=1}^\infty \dfrac{1}{(n-1)!}\ x^{n-1} y^{(n-1,h)}\bigl(y-(n-1)h\bigr)\\
&=y\sum_{n=1}^\infty \dfrac{1}{(n-1)!}\ x^{n-1} y^{(n-1,h)}
-h\sum_{n=2}^\infty \dfrac{1}{(n-2)!}\ x^{n-1} y^{(n-1,h)}\\
&=y\ e_h(x,y)-hx\frac{\partial}{\partial x}e_h(x,y)\ ,
\endaligned
$$
wherefrom we get
$$
(1+hx)\frac{\partial}{\partial x}\ e_h(x,y)=y\ e_h(x,y),
$$
i.e.,
$$
\frac{d}{d_h x}\ e_h(x,y)=y\ e_h(x,y).\quad \Box
$$

\begin{thm}
The function $x\mapsto \exp_h(x,y)$ is the eigenfunction of operator $\dfrac d{d^h x}$ with eigenvalue $y$.
\end{thm}

\noindent{\it Proof.} If we apply differential operator $\dfrac{\partial}{\partial x}$ on func\-tion $\exp_h(x,y)$, we obtain
$$
\aligned
\frac{\partial}{\partial x}\exp_h(x,y)
&=\frac{y}{h}\bigl(hx+\sqrt{1+h^2 x^2}\bigr)^{\frac{y}{h}-1}\Bigl(h+\frac{h^2 x}{\sqrt{1+h^2 x^2}}\Bigr)\\
&=\frac{y}{\sqrt{1+h^2 x^2}}\bigl(hx+\sqrt{1+h^2 x^2}\bigr)^{\frac{y}{h}}\\
&=\frac{y}{\sqrt{1+h^2 x^2}}\ \exp_h(x,y),
\endaligned
$$
i.e.,
$$
\Bigl(\sqrt{1+h^2 x^2}\ \frac{\partial}{\partial
x}\Bigr)\exp_h(x,y)=y\exp_h(x,y).
$$
Hence,
$$
\frac{d}{d^h x}\ \exp_h(x,y)=y\ \exp_h(x,y).\quad \Box
$$

Finally, let us consider behavior of deformed exponential functions related to differentiating over the second variable. In certain sense,
we can conclude that they are "deformed" eigenfunctions of operator $\partial/\partial y$. Really, the following is valid:
$$
\aligned
\frac {\partial}{\partial y}e_h(x,y)&=\{x\}_h\ e_h(x,y),\\
\frac {\partial}{\partial y}\exp_h(x,y)&=\{x\}^h\ \exp_h(x,y).
\endaligned
$$

\medskip

\noindent {\bf Acknowledgement.} Supported by Ministry of Sciences and Technology of Republic Serbia, the projects  No $144023$ and No
$144013$.

\end{document}